\documentclass[12pt]{article}

\usepackage{latexsym}

\begin{document}

\thispagestyle{empty}

\vskip 20pt
\begin{center}
{\Large Some New Exact van der Waerden numbers} \vskip 15pt
{\bf  Bruce Landman}\\
{\it Department of Mathematics}\\
{\it  University of West Georgia,
Carrollton, GA 30118}\\{\tt landman@westga.edu}\\
\vskip 20pt
{\bf Aaron Robertson} \\
{\it Department of Mathematics}\\
{\it Colgate University, Hamilton, NY}\\ {\tt aaron@math.colgate.edu}\\
\vskip 20pt
{\bf Clay Culver} \\
{\it University of West Georgia, Carrollton, GA 30118} \vskip 10pt
\end{center}

\vskip 30pt

\begin{abstract}{\footnotesize
For positive integers $r,k_0,k_1,...,k_{r-1},$ the van der Waerden number
$w(k_0,k_1,...,k_{r-1})$ is the least positive integer $n$ such that whenever $\{1,2,\dots,n\}$ is
partitioned into $r$ sets
$S_{0},S_{1},...,S_{r-1}$, there is some $i$ so that $S_i$ contains a $k_i$-term arithmetic progression. We
find several new exact values of $w(k_0,k_1,...,k_{r-1})$. In addition,
for the situation in which only one value of
$k_i$ differs from 2, we give a precise formula for the van der Waerden function (provided this one value
of
$k_i$ is not too small)}.
\end{abstract}

\vskip 20pt

\section{Introduction}

 A well-known theorem of van der Waerden [9] states that for any positive integers $k$ and
$r$, there exists a least positive integer, $w(k;r)$, such that any $r$-coloring of
$[1,w(k;r)]=\{1,2,\dots,w(k;r)\}$ must contain a monochromatic $k$-term arithmetic progression
$\{x,x+d,x+2d,\dots,x+(k-1)d\}$. Some equivalent forms of van der Waerden's theorem may be found in
 in [7].
It is easy to see that the existence of the van der Waerden numbers $w(k;r)$ implies the existence of the (more
general) van der Waerden numbers $w(k_0,k_1,...,k_{r-1};r)$ which are defined as follows.

\noindent {\bf Definition}. {\em Let $k_0,k_1,...,k_{r-1}$ be positive integers. The van der Waerden number
$w(k_0,k_1,...,k_{r-1};r)$ is the least positive integer $n$ such that every $r$-coloring $\chi:[1,n]
\rightarrow \{0,1,...,r-1\}$ admits, for some $i$, $0 \leq i \leq r-1$, a $k_i$-term arithmetic progression of
color $i$.}

If the value of $r$ is clear from the context, we will denote the van der Waerden number $w(k_0,k_1,...,k_{r-1};r)$ more simply
by $w(k_0,k_1,...,k_{r-1})$.

For example, $w(4,4)= w(4,4;2)$ has the same meaning as the classical van der Waerden number $w(4;2)$;
$w(3,3,3,3)$ has the same meaning as $w(3;4)$; and $w(5,4,7;3)=w(5,4,7)$
represents the least positive $n$ such that every
(red,blue,green)-coloring of $[1,n]$ yields a red 5-term arithmetic progression, a blue 4-term arithmetic
progression, or a green 7-term arithmetic progression.

The study of these ``mixed" van der Waerden numbers apparently has received relatively little attention,
especially when compared to, say, the classical (mixed) graph-theoretical Ramsey numbers $R(k_1,k_2,...,k_r)$.
It is easy to calculate by hand that $w(3,3)=9$. Other non-trivial values of the van der Waerden numbers were
published by Chv\'{a}tal [4], Brown [3], Stevens and Shantaram [8], Beeler and O'Neil [2], and Beeler [1], in
1970, 1974, 1978, 1979, and 1983, respectively.

The purpose of this note is to expand the table of known van der Waerden numbers. In Section 2, we
present several new van der Waerden numbers. In Section 3, we give a
formula for the van der Waerden numbers of the type
$w(k,2,2,...,2;r)$, i.e., where all but one of the $k_i$'s equal 2, provided $k$ is large enough in relation
to
$r$.

\section{Some New Values}

Table 1 below gives all known non-trivial van der Waerden numbers, $w(k_0,k_1,...,k_{r-1})$, where at least two
of the $k_i$'s exceed two (the cases where only one $k_i$ exceeds two are handled, more generally, in the next
section). The entries in the table due to Chv\'atal, Brown, Stevens and Shantaram, Beeler and O'Neil, or Beeler,
are footnoted by a, b, c, d, or e, respectively. The entries that we are presenting here as previously unknown
values are marked with the symbol *.

To calculate most of the new van der Waerden numbers, we used a slightly modified version of the ``culprit"
algorithm, introduced in [2]. In certain instances, we used a very simple backtracking algorithm, which can be
found in [8]. These, and a third algorithm, are also described in [7].

\begin{center}

\begin{tabular}{|c|c|c|c|c|c|c|} \hline
$r$&$k_{0}$&$k_{1}$&$k_{2}$&$k_{3}$&$k_{4}$&$w$\\ \hline \hline $2$&$3$&$3$&-&-&-&$9$ \\ \hline
$2$&$4$&$3$&-&-&-&$18^{a}$\\
\hline
$2$&$4$&$4$&-&-&-&$35^{a}$ \\ \hline $2$&$5$&$3$&-&-&-&$22^{a}$\\ \hline $2$&$5$&$4$&-&-&-&$55^{a}$ \\
\hline
$2$&$5$&$5$&-&-&-&$178^{c}$ \\
\hline $2$&$6$&$3$&-&-&-&$32^{a}$ \\ \hline $2$&$6$&$4$&-&-&-&$73^{d}$ \\ \hline
$2$&$7$&$3$&-&-&-&$46^{a}$\\
\hline $2$&$7$&$4$&-&-&-&$109^{e}$
\\ \hline $2$&$8$&$3$&-&-&-&$58^{d}$ \\ \hline $2$&$9$&$3$&-&-&-&$77^{d}$ \\ \hline
$2$&$10$&$3$&-&-&-&$97^{d}$ \\
\hline $2$&$11$&$3$&-&-&-&$114^{*}$ \\ \hline $2$&$12$&$3$&-&-&-&$135^{*}$ \\ \hline
$2$&$13$&$3$&-&-&-&$160^{*}$ \\ \hline $3$&$3$&$3$&$2$&-&-&$14^{b}$
\\ \hline $3$&$3$&$3$&$3$&-&-&$27^{a}$ \\ \hline $3$&$4$&$3$&$2$&-&-&$21^{b}$
\\
\hline $3$&$4$&$3$&$3$&-&-&$51^{d}$ \\ \hline $3$&$4$&$4$&$2$&-&-&$40^{b}$ \\ \hline
$3$&$4$&$4$&$3$&-&-&$89^{*}$
\\
\hline $3$&$5$&$3$&$2$&-&-&$32^{b}$\\ \hline $3$&$5$&$3$&$3$&-&-&$80^{*}$ \\ \hline $3$&$5$&$4$&$2$&-&-&$71^{b}$
\\
\hline $3$&$6$&$3$&$2$&-&-&$40^{b}$ \\ \hline $3$&$6$&$4$&$2$&-&-&$83^{*}$ \\ \hline
$3$&$7$&$3$&$2$&-&-&$55^{*}$
\\
\hline $4$&$3$&$3$&$2$&$2$&-&$17^{b}$ \\ \hline $4$&$3$&$3$&$3$&$2$&-&$40^{b}$ \\ \hline
$4$&$3$&$3$&$3$&$3$&-&$76^{d}$
\\
\hline $4$&$4$&$3$&$2$&$2$&-&$25^{b}$ \\ \hline $4$&$4$&$3$&$3$&$2$&-&$60^{*}$ \\ \hline
$4$&$4$&$4$&$2$&$2$&-&$53^{b}$
 \\ \hline $4$&$5$&$3$&$2$&$2$&-&$43^{b}$ \\
\hline
$4$&$6$&$3$&$2$&$2$&-&$48^{*}$ \\
\hline
$4$&$7$&$3$&$2$&$2$&-&$65^{*}$ \\
\hline
$5$&3&3&2&2&2&$20^*$\\ \hline
$5$&3&3&3&2&2&$41^*$\\ \hline
\end{tabular}

{\bf TABLE 1. Van der Waerden Numbers}
\end{center}

An $r$-coloring of an interval $[a,b]$ is said to be
$(k_0,k_1,...,k_{r-1})$-{\em valid} (or simply {\em valid} if
the specific  $r$-tuple is clear from the context) if, for each $i=0,1,...,r-1$, the coloring avoids $k_i$-term
arithmetic progressions of color $i$.

For each of  the new van der Waerden numbers $w= w(k_0,k_1,...,k_{r-1})$ we computed, the program also outputted
all maximal length $(k_0,k_1,...,k_{r-1})$-valid $r$-colorings (i.e., all valid colorings of length $w-1$). We
now describe all of these maximal-length valid colorings. For convenience, we will denote a coloring as a string
of colors. For example, the coloring $\chi$ of $[1,5]$ such that $\chi(1)=\chi(2)=0$ and
$\chi(3)=\chi(4)=\chi(5)=1$, will be denoted by the string 00111. Further, a string of $t\geq 2$ consecutive
$i$'s will be denoted by $i^t$. Thus, for example, $0^{2}1^{3}$ represents 00111. In counting the number of
valid colorings, we will consider two colorings to be the same if one can be obtained
from the other by a
renaming of the colors. For example, we consider 11001 and 00110 to be the same coloring of [1,5].

\noindent $\mathbf{w(11,3)}$:
 Corresponding to $w(11,3)=114$, there are thirty different valid
2-colorings of [1,113]. These are:
\renewcommand{\arraystretch}{1.25}
$$
\begin{array}{ll}
\,\, 0^{10}10^{10}1^{2}0^{6}10^{10}1^{2}01^{2}0^{9}
1010^{6}10^{9}1^{2}0^{2}10^{10}10^{6}1^{2}0^{4}10^{7}a0^{2}b,\\
\mbox{where at least one of $a$ and $b$ equals 1;}
\\
\,\, 0^{7}c010^{7}1^{2}0^{3}10^{7}1^{2}0^{4}10^{10}1010^{6}10^{9}
10^{2}10^{3}1010^{2}d0^{7}1010^{3}1^{2}0^{10}1e0^{6}1,\\
\mbox{where the only restriction on $c$, $d$, and $e$ is that $c$ and $e$ are
not both 1;}
\\
\,\, 0^{9}10^{3}1010^{2}f010^{8}10^{3}g1010^{10}10^{8}10^{9}1^{2}0^{5}
1010^{10}101^{2}010^{10}1010^{3}h0^{4},\\
\mbox{where the only restriction on $f$, $g$, and $h$ is that $f$ and $g$ are
not both 0;}
\end{array}
$$
and the fifteen colorings obtained by reversal of the
fifteen colorings described above.

\noindent $\mathbf{w(12,3)}$: Corresponding to $w(12,3)=135$, there is only one valid 2-coloring of $[1,134]$:
$$0^{9}10^{8}10^{9}
10^{2}10^{3}1010^{7}10^{2}1010^{3}10^{11}1^{2}0^{11}10^{3}1010^{2}10^{7}1010^{3}10^{2}10^{9}10^{8}10^{9}$$
(it is its own reversal).

\noindent $\mathbf{w(13,3)}$: There are twenty-four different (13,3)-valid 2-colorings of [1,159]:
$$
\begin{array}{l}
\,\,0^{12}10^{2}a0^{7}101^{2}0^{12}10^{2}1^{2}0^{6}10^{3}10^{9}10^{2}1^{2}
0^{4}10^{2}10^{4}10^{8}10^{10}10^{4}10^{7}1^{2}0^{2}10^{9}10^{3}1^{2}0^{9}10^{12}10^{2},
\\
\mbox{where $a$ may be 0 or 1;}
\\
\,\,0^{2}10^{4}10^{5}10b0^{10}10^{5}1^{2}010^{10}10^{7}10^{3}10^{2}
10^{9}10^{2}c0^{3}1010^{11}10^{4}10^{2}10^{10}10^{6}1^{2}0^{11}1^{2}0^{4}10^{12}10^{6}10,\\
\mbox{with no
restrictions on $b$ or $c$;}
\\
\,\,0^{2}1010^{4}d0^{5}10^{7}10^{6}10^{3}
0^{5}1^{2}010^{5}10^{10}e010^{6}10^{2}10^{12}101^{2}0^{3}10^{7}10^{5}10^{4}10^{5}10
10^{6}10^{12}0^{5}10^{2}10^{7}f01,\\
\mbox{where the only restriction is that $d$ and $e$ are not both 1;}
\end{array}
$$
and the twelve others obtained from these twelve via reversal.

\noindent $\mathbf{w(4,4,3)}$: There are four (4,4,3)-valid 3-colorings of [1,88]:
$$
0^2a101^30^31^302121^22^2010^2202^201^22^20^3101^30^31^302121^32010^2202^20102^20^21^201^30^3
21^201^20
$$
where $a$ may be $0$ or $1$, along with the two reversals of these colorings.

\noindent $\mathbf{w(5,3,3)}$: There are forty-two different $(5,3,3)$-valid 3-colorings of [1,79]:

\noindent
Twenty-one are of the form
$$0^{2}a102b0^{2}101^{2}210^{3}10^{4}20^{3}212^{2}020^{2}1^{2}021^{2}02^{2}0^{2}101^{2}210^{3}10^{4}20^{3}212^{2}020^{2}c102d0^{2},$$
 where $ab \neq 00$, $cd \neq 00$, $b \neq 1$, $c \neq 2$, and $ad \neq
11$; by reversal of these colorings,
the other twenty-one colorings are obtained (note: of course, since $k_{2}=k_{3}$, interchanging of the colors
1 and 2 will not change the validity of any coloring, so we are not
considering the forty-two valid colorings of
[1,79] obtained in this manner, from those already listed, to be
additional colorings).

\noindent $\mathbf{w(6,4,2)}$: There are twelve different $(6,4,2)$-valid 3-colorings of [1,82]:
$$
\begin{array}{l}
0^{2}10101^{2}0^{3}10101^{2}010^{3}1^{2}0^{3}10^{5}1^{2}01^{2}21^{2}0^{4}10^{2}10^{3}1^{2}010^{2}10^{3}1
0101^{2}0^{5}10^{3}1^{2}0^{2};
\\

0^{4}10^{5}1^{2}01010^{4}101^{2}0^{3}10^{4}1^{3}010^{4}10^{3}101^{2}21^{2}010^{4}10^{2}1a0^{3}
1^{3}010^{5}1b01c0,\\
\mbox{where either $a=0$ or else $(a,b,c)=(1,0,0)$ (this gives five colorings);}
\end{array}
$$
\noindent and six more by reversal of the above  colorings.

\noindent $\mathbf{w(7,3,2)}$: To describe the  $(7,3,2)$-valid 3-colorings of $[1,54]$, we adopt the following
notation. For $1 \leq i \leq 10$, denote by $S_i$ the coloring
$$0^{6}1^{2}0^{3}10^{5}101^{2}0^{3}1010^{6}20^{4}a010^{2}1^{2}0^{2}b0cd0e, $$
where $(a,b,c,d,e)$ is the $i$th element of the $10$-element set
$\{(0,0,0,1,1), (0,0,1,0,1),$ $
 (0,0,1,1,1),  (0,1,0,0,0), (0,1,0,0,1),$ $(0,1,1,0,0), (0,1,1,0,1),$
$(1,1,0,0,0),  (1,1,1,0,0),$ $(0,1,0,1,1)\}.$ Then the valid colorings of
$[1,54]$ are: $$ \begin{array}{l}
\,\,  0S_{i}, \hskip 10pt 1 \leq i \leq 9; \\
\,\,  1S_{i}; \hskip 10pt \mbox{where $i \in \{1,4,5,8,10\}$}; \\
\,\, 0^{3}10^{5}1^{2}0^{2}10^{4}101^{2}0^{5}20^{6}1010^{3}1^{2}010^{5}10^{3}a, \mbox{where $a$ may be 0 or 1;}
\\
\,\,0^{4}10^{5}10^{0}10^{4}101^{2}0^{5}20^{5}1^{2}0^{5}1^{2}01^{2}0^{6}10^{2};\\
\,\, 010^{3}0^{4}10^{2}10^{4}10^{2}20^{6}101^{2}0^{3}1^{2}0^{6}10^{3}1;
\end{array}
$$
\noindent
 and eighteen obtained by reversal of those already listed, giving a total of thirty-six valid
colorings.

\noindent $\mathbf{w(4,3,3,2)}$: There are eight different $(4,3,3,2)$-valid 4-colorings of [1,59]:
$$1012^{2}01020^{2}a01^{2}2021^{2}0^{3}320^{3}20^{2}1012^{2}0121^{2}0120^{2}2^{2}0^{3}
b021^{2},$$ \noindent where $a,b \in \{1,2\}$, and four others obtained by reversal.

\noindent
$\mathbf{w(6,3,2,2)}$:
There are twenty-eight different $(6,3,2,2;4)$-valid colorings of [1,47]:
$$
\begin{array}{l}
\,\, 00001000211011000011010000300001000110000010011
\\
\,\,0000101001000001021100000100100031001000a011010, \,\,a \in \{0,1\}
\\
\,\, 0000101001000001021100100000100131001000001001a, \,\,a \in \{0,1\}
\\
\,\,00010000112110000110100003000010001100000100110
\\
\,\,0100001001000001121100000100100031001000001101a, \,\,a \in \{0,1\}
\\
\,\,0100101001000001021100000100100031001000a011010, \,\,a \in \{0,1\}
\\
\,\,0100101001000001021100100000100131001000001001a, \,\,a \in \{0,1\}
\\
\,\,1a001010010000210001001000001301010010000010011, \,\,a \in \{0,1\}
\end{array}
$$
\noindent plus fourteen more colorings obtained by reversing the above fourteen colorings.

\noindent $\mathbf{w(7,3,2,2)}$: There are five different $(7,3,2,2)$-valid colorings of [1,64]:
$$0^310^61^20a10^610^210^512310^61^201^20^610^210^510^3$$
\noindent where $a \in \{0,1\}$, along with the reversals of these; and
$$0^210^310^61^201^20^220^61010^41010^630^21^201^20^610^310^2,$$
\noindent
which is its own reversal.

\noindent
$\mathbf{w(3,3,2,2,2)}$:  There are five different $(3,3,2,2,2;5)$-valid
colorings of [1,19]:
$$
\begin{array}{l}
\,\,0^210^21^2231410^21^20^21,\,\, \mbox{and its reversal}; \\
\,\,0^21^20^21^22340^21^20^21^2;\\
\,\,0^21^221^20^230^21^241^20^2;\\
\,\,0101^20102340101^2010.
\end{array}
$$

\noindent $\mathbf{w(3,3,3,2,2)}$: There are forty-two different $(3,3,3,2,2)$-valid 5-colorings of [1,40]:
$$
\begin{array}{l}
\,\,0xy^22y2x^22^2x^22y2y^2x3xy^22y2x^22^2x^22y2y^2xb4, \\
 \mbox{where $xy=01$ or $xy=10$, and $b \in \{0,1,2\}$;} \\
\,\,010x^2yxy0^2y^20^2yxyx^2ab0x^2yxy0^2y^20^2yxyx^2cd, \\
 \mbox{where $xy=12$ or $xy=21$,
$abc=  034$ or $abc=340$, and $d \in \{ 0,1,2\}$;}\\
\,\, 012^20201^20^21^20202^21312^20201^20^21^20202^21a4, \mbox{where $a\in\{0,1,2\}$;}
  \end{array}
$$
along with the reversals of the (twenty-one) colorings above.

\section{When All But One $k_{i}$ Equals 2}

Consider $w(k,2,2,...,2;r)$. For ease of notation, we will denote this function more simply as $w_{2}(k;r)$. We
give an explicit formula for $w_{2}(k;r)$ provided $k$ is large enough in relation to $r$.

We adopt the following notation. Let $p_{1} < p_{2} < \cdots $ be the sequence of primes. For $r \geq 2$, 
denote by $\pi(r)$ the number of primes not exceeding $r$, and denote by $\#r$ the product $p_{1}p_{2}\cdots
p_{\pi(r)}$. For $k,r \geq 2$, let $j_{k,r} = \min\{j \geq 0: \gcd(k-j,\#r)=1\}$, and $\ell_{k,r}=\min\{\ell
\geq 0: \gcd(k-\ell,\#r)=r\}$.

\renewcommand{\arraystretch}{1.0}

\noindent {\bf Theorem} Let $k > r \geq 2$. Let $j=j_{k,r}$,  $\ell=\ell_{k,r}$, and $m=\min\{j,\ell\}$.
\begin{enumerate}
\item[I.] $w_{2}(k;r) = rk$ if $j=0$. \item[II.] $w_{2}(k;r) = rk-r+1$ if either (i) $j=1$; or (ii) $r$ is prime
and $\ell=0$. \item[III.] If $r$ is composite and $j \geq 2$, then $w_2(k;r) \geq rk-j(r-2)$, with equality
provided either (i) $j=2$ and $k \geq 2r-3$; or (ii) $j \geq 3$ and $k \geq (\pi(r))^{3}(r-2)$.
 \item[IV.] If
$r$ is prime, $j\geq 2$, and  $\ell \geq 1$, then $w_2(k;r) \geq rk-m(r-2)$, with equality provided
 either (i) $\ell =1$, (ii) $m=2$ and $k \geq 2r-3$, or (iii)
$m \geq 3$ and $k  \geq (\pi(r))^{3}(r-2)$.
\end{enumerate}

 \noindent {\em Proof}. (I) To show $w_{2}(k;r) \leq rk$, let $g: [1,rk] \rightarrow \{0,1,2,...,r-1\}$ be any
$r$-coloring. If $g$ is valid, then no $k$ consecutive elements of $[1,rk]$ have color 0, and for each $i$, $1
\leq i \leq r-1$, there is at most one element with color $i$. This is not possible, since these last two
conditions would imply that no more than $r-1+r(k-1) =rk-1$ integers have been colored. Hence, $g$ is not valid.

To show that $rk$ is also a lower bound, we show that  the following coloring of $[1,rk-1]$: is valid:
$$0^{k-1}10^{k-1}20^{k-1}3 \ldots 0^{k-1}(r-1)0^{k-1}.$$ Obviously, this coloring admits no 2-term monochromatic
arithmetic progression having color other than 0. For a contradiction, assume there is a $k$-term arithmetic
progression $A=\{a+id: 0 \leq i \leq k-1\}$ of color 0. Since $k > r$, we have $1 < d \leq r$. Hence, by
hypothesis, $(d,k)=1$, so $A$ is a complete residue system modulo $k$. Therefore there is some $x \in A$ with $x
\equiv 0$(mod $k$). Then $x$ does not have color 0, a contradiction.

 (II) (i) We first note that $k$ must be even, for otherwise $(k-1,\#r) \neq
1$.
 We prove that $r(k-1)+1$ is an upper bound by induction on $r$. First, let
$r=2$.
 Let $g$ be an arbitrary (0,1)-coloring of $[1,2k-1]$. If $g$ is $(k,2)$-valid and
if
 there is at most one integer having color 1, then $g$ must be the coloring
$0^{k-1}10^{k-1}$
 (since there cannot be $k$ consecutive integers with color 0). This contradicts the validity of $g$ since $g(\{2i-1: 1 \leq
 i \leq k\})=0$. This establishes the desired upper bound for $r=2$.

 Now let $r \geq 3$ and assume that for all $s$, $2 \leq s < r$, we have $w_{2}(k;s) \leq
s(k-1)+1$ whenever $j_{k,s}=1$. Let $h$ be any
 $r$-coloring of $[1,r(k-1)+1]$, with $j_{k,r}=1$. Without loss of generality,
we may assume $h = B_{1}1B_{2}2B_{3}3...B_{r-1}(r-1)B_{r}$, where each $B_{i}$ is a (possibly empty) string of
0's. For each $i$, $1 \leq i \leq r-1$, let $x_i$ be the integer for which $h(x_i)=i$.

For a contradiction, assume $h$ is $(k,2,2,...,2;r)$-valid.  Note that $|B_i| \leq k-1$ for $1 \leq i \leq r-2$,
so that $|B_r| \geq 1$.
 Now, for each $t$, $1 \leq t
\leq r-1$, since all primes $p \leq t$ do not divide $k-1$ and since $k$ is even,  $j_{k,t}=1$. By the
induction hypothesis, $ w_{2}(k,t) \leq t(k-1)+1$, and therefore we may
assume that $h$ admits at least $t+1$ colors
within $[1,t(k-1)+1]$ (or else $h$ would not be valid). Thus, $x_{t} \leq t(k-1)+1$ for $1 \leq t \leq r-1$.
By a symmetric argument (considering intervals of the form
$[(r-t)(k-1)+1,r(k-1)+1]$), it follows that $x_{t} \geq
t(k-1)+1$ for all $t$, $1 \leq t \leq r-1$. Hence for each $t$, $1 \leq t \leq r-1$, we
have
 \[
 x_{t}=t(k-1)+1.\]

Consider $S=\{1+ir: 0 \leq i \leq (k-1)\}$. We claim that $S$ is monochromatic of color 0, the truth of which
contradicts the assumption that $h$  is $(k,r)$-valid. Let $1+cr \in S$. If $h(1+cr) \neq 0$, then $cr = n(k-1)$
for some $n$, $1 \leq n \leq r-1$. However, $(k-1)$ cannot divide $cr$
 since $\gcd (k-1,r)=1$ and $c
  \leq k-2$ (since $|B_{r}| \geq 1$). Hence, $h(S)=0$, which completes the proof that
$r(k-1)+1$
  is an upper bound.

The coloring $0^{k-1}10^{k-2}20^{k-2} \ldots (r-2)0^{k-2}(r-1)0^{k-2}$ is $(k,2,2,...,2;r)$-valid, which
establishes that $r(k-1)+1$ is also a lower bound. To see that this coloring is valid, note that if $S=\{a+id: 0
\leq i \leq k-1\} \subseteq [1,r(k-1)]$ is an arithmetic progression, then since $d \leq r-1$, we have
$\gcd(k-1,d)=1$. Thus, $S$ represents all residue classes modulo $k-1$, in particular the class $1$(mod
$(k-1)$), with $a$ and $a+(k-1)d$ in the same class. Hence, since no member of $\{1+i(k-1): 1 \leq i \leq r-1\}$
has color 0, $S$ cannot be monochromatic.

\noindent (II) To show $r(k-1)+1$ is a lower bound, we show that the coloring $$ 0^{k-1} 1 0^{k-1} 2 ...(r-2)
0^{k-1} (r-1) 0^{k-r}$$ is valid. Let $A$ be a $k$-term arithmetic progression. Obviously, there is no $k$-term
arithmetic progression having color 0 and gap $d=1$. So we may assume $2 \leq d < r$. Since $A$ forms a complete
residue system modulo $k$, some member of $A$ is a multiple of $k$. Hence $A$ is not monochromatic, so the
coloring is valid.

To complete the proof, let $f$ be any $r$-coloring of $[1,r(k-1)+1]$, and assume (for a contradiction) that $f$
is valid. By the proof of Part I, we may assume that for $1 \leq s \leq r-1$, the interval $[1,ks]$ assumes at
least $s+1$ colors, and the interval $[k-r+2,r(k-1)+1]$ assumes all $r$ colors. Since all non-zero colors occur
at most once, $[ks+1,r(k-1)+1]$ assumes at most $r-s$ colors for $1 \leq s \leq r-1$. Similarly,
$[1,k(r-s)-r+1]$ assumes at most $r-s$ colors. From this we may conclude that $f([sk+1,ks+k-r+1])=0$ for each
$s$, $0 \leq s \leq r-1$. Since $r|k$, we obtain $f(1+jr)=0$ for $j=0,1,...,k-1$, a contradiction.

\noindent(III) Note that the coloring
$$
0^{k-1}10^{k-j-1}20^{k-j-1}3\dots0^{k-j-1}(r-1)0^{k-1} $$ is a valid coloring of $[1,n]=[1,kr-j(r-2)+1]$. (The
argument is very similar to that used to establish the lower bounds in Parts I and II of the theorem.)

To complete the proof,
we will show that, subject to the stated restrictions on $k$,
the above coloring is the {\em only}
valid $r$-coloring of $[1,n]$.
The desired result then follows since the coloring clearly cannot be extended to a
valid $r$-coloring of $[1,n+1]$.

Let $\tau: [1,n] \rightarrow \{0,1,\dots,r-1\}$ be valid.  We may assume that
$$
\tau = 0^{\alpha_1}10^{\alpha_2}2\dots 0^{\alpha_{r-1}}(r-1) 0^{\alpha_r},
$$
\noindent with $0 \leq \alpha_i\leq k-1$ for $i = 1,2,\dots,r-1$. For $i=1,2,\dots,r-1$, let $\tau(y_i)=i$. Let
$y_0=0$ and $y_r=n+1$. Define
$$
B_i=\{x:y_{i-1}<x<y_i\}
$$
\noindent
for $i=1,2,\dots,r$, so that $|B_i|=\alpha_i$.

Let $b \geq 1$, $t \geq 2$, and assume that $b+t \leq r$.  We show that there do not exist $i_1$ and $i_2$ such
that
\begin{equation}
b \leq i_1 <i_2 \leq b+t-1 \mbox{ and } y_{i_1} \equiv y_{i_2} \, (\mbox{mod }t).
\end{equation}

Let
$$
a= \sum_{i=b}^{b+t} \alpha_i.
$$
We know that $\sum_{i=1}^r \alpha_i = kr - j(r-2)-1-(r-1) = r(k-1)-j(r-2)$.
Since $\alpha_i \leq k-1$ for all $i \in [1,r]$, we have
\begin{equation}
a \leq r(k-1)-j(r-2)-(k-1)(r-t-1) = (k-1)(t+1)-j(r-2).
\end{equation}

Now consider the interval
$$
I = \left(\bigcup_{i=b}^{b+t} B_i \right)\bigcup \{y_b,y_{b+1},\dots,y_{b+t-1}\}.
$$
We claim that $|I| \geq kt$ (for both Cases (i) and (ii)). Using (2), we see that $|I| \geq (k-1)(t+1)-2(r-2)+t
= kt+k-1-2(r-2)$. Since $k\geq 2r-3$, the claim holds in Case (i) .

For Case (ii) we turn to a result of Jurkat and Richert [6].  Let $C(r)$ be the length of the longest string of
consecutive integers, each divisible by one of the first $\pi(r)$ primes (this is a particular case of what is
known as Jacobsthal's function). Jurkat and Richert showed that $C(r) < (\pi(r))^2 \exp(\log(\pi(r))^{13/14})$.
(As an aside, better asymptotic bounds have been proved: Iwaniec [5] showed that $C(r) \ll \pi(r)^2
\log^2(\pi(r))$.) Using the slightly weaker bound $C(r) < \pi(r)^3$, we have that $j<\pi(r)^3$. Hence, from (2)
we see that $|I| > kt+k-1-\pi(r)^3(r-2)$.  Since $k\geq\pi(r)^3(r-2)$, we have $|I| \geq kt$.

Having established that $|I| \geq kt$, now assume, for a contradiction, that there exist $i_1$ and $i_2$
satisfying (1). Then there is some $c$, $0 \leq c \leq t-1$, so that no member of
$\{y_m,y_{m+1},\dots,y_{m+t-1}\}$ is congruent to $c$ modulo $t$.  Since $\{y_m,y_{m+1},\dots,y_{m+t-1}\}$ are
the only members of $I$ {\it not} of color $0$, the set $\{x \in I: x \equiv c \,(\mbox{mod }t)\}$ is an
arithmetic progression of length at least $k$ and of color $0$, a contradiction.  Hence, (1) is true.

By (1) we have, in particular, that for all $i \in [1,r-2]$, $y_{i+1}-y_i \not \equiv 0 \,(\mbox{mod }t)$ for
each $t = 2,3,\dots,r-1$.  By assumption, for all $u \in [k-j+1,k]$, there exists $t \in [2,r-1]$ that divides
$u$.  Hence, for each $i \in [1,r-2]$ we must have $y_{i+1}-y_i \leq k-j$.  Thus,
$$
\alpha_i \leq k-j-1 \hskip 10pt \mbox{for}\,\, 2 \leq i \leq r-1.
$$
Since
$$
\sum_{i=2}^{r-1} \alpha_i \geq (r-2)(k-j-1)
$$
we must have $\alpha_i = k-j-1$ for $i \in [2,r-1]$, and hence $\alpha_1=\alpha_r=k-1$.

 \noindent (IV) The proof is essentially the same as that for Part III, where
we use the $r$-coloring $$ 0^{k-1}10^{k-m-1}20^{k-m-1}3\dots0^{k-m-1}(r-1)0^{k-1}, $$ which is a valid coloring
of $[1,n]=[1,kr-m(r-2)+1]$. For this reason, we include the details  for only Case (i),  where we   need only
the restriction that $k > r$.

To prove the result for Case (i), note first that since $\gcd(k-1,\#(r-1))=1$, by  II.ii of this theorem, we
know that $w_{2}(k;r-1)= k(r-1)-r+2$. Using the notation from the proof of Part III , we thus may assume that
$y_1 \geq k$ and $y_{r-1} \leq k(r-1)-r+2$. This implies that $$ \sum_{i=2}^{r-1} \alpha_{i} = (r-2)(k-2),$$ and
therefore, for any $u$, $1 \leq u \leq r-2$, we have
$$\sum_{i=u}^{u+2}\alpha_i \geq 3(k-1)-r+2.$$ Now if some $\alpha_i$ is odd, where $2 \leq i \leq r-1$, then
there is a 0-colored arithmetic progression with gap 2 and length at least \[\frac{3k-r-1}{2} + \frac{1}{2} =
\frac{3k-r}{2} \geq k.\] Otherwise, if each $\alpha_{i}$, $2 \leq i \leq r-1$, is even, then since $k$ is even
(by hypothesis) and $\alpha_i \leq k-1$, we must have $\alpha_i =k-2$ for $i=2,\dots,r-1$.
 \hfill
  $\Box$

\noindent {\bf Remarks}. It is easy to see that the converse of Part I is also true.
It is also clear that, for any fixed $r$, the values of
$w_2(k;r)$ given by the theorem may be expressed in terms of the
residue classes of $k$ modulo $\#r$.  This is convenient for small values
of
$r$.  As examples, for $k \geq 4$, we have $w_2(k;3) =  3k$ if $k \equiv
\pm 1($mod $6), w_2(k;3) = 3k-1$ if $k \equiv 3($mod $6)$, and
$w_2(k;3)=3k$ otherwise; and, for $k \geq 5$, we have $w_2(k;4)=4k-3$ if
$k \equiv 0$ or $2 ($ mod $6), w_2(k;4) = 4k-4$ if $k \equiv 3($mod $6),
w_2(k;4)= 4k-6$ if $k \equiv 4($mod $6)$ (the theorem gives this last
equality only for $k \geq 16$, but it is easy to show that
it holds provided $k \geq 5$), and $w_2(k;4) = 4k$ otherwise.
We believe that the restrictions on the magnitude of $k$ (in
relation to
$r$) in Parts III.ii and IV.iii can be weakened substantially. We also
would like to obtain precise formulas in the other cases, but with weaker 
restrictions on the magnitude of $k$.

\section*{Acknowledgement}  We thank Carl Pomerance for
pointing us to some relevant literature.

\section*{References}
\parindent=0pt \footnotesize

[1] M. Beeler, A new van der Waerden number, {\em Discrete Applied Math.} {\bf 6} (1983), 207.

[2] M. Beeler and P. O'Neil, Some new van der Waerden numbers, {\em Discrete Math.} {\bf 28} (1979), 135-146.

[3] T.C. Brown, Some new van der Waerden numbers (preliminary report), {\em Notices American Math. Society} {\bf
21} (1974), A-432.

[4] V. Chv\'atal, Some unknown van der Waerden numbers, in {\em Combinatorial Structures and their
Applications}, 31-33, Gordon and Breach, New York, 1970.

 [5] H. Iwaniec, On the error term in the linear
sieve, {\it Acta Arith.} {\bf 19} (1971), 1-30.

[6] W. B. Jurkat and H. E. Richert, An improvement of Selberg's sieve method I, {\it Acta Arith.} {\bf 11}
(1965), 217-240.

[7] B. Landman and A. Robertson, Ramsey Theory on the Integers, {\it American Math. Society}, Providence, RI,
317pp., 2004.

[8] R. Stevens and R. Shantaram, Computer-generated van der Waerden partitions, {\em Math. Computation} {\bf 32}
(1978), 635-636.

 [9] B.L. van der Waerden, Beweis einer Baudetschen Vermutung, {\em Nieuw Arch. Wisk.} {\bf 15} (1927), 212-216.

\end{document}